\documentclass[12pt]{amsart}
\numberwithin{equation}{section}

\usepackage{amssymb}
\def\ca{{\mathcal A}}

\def\cd{{\mathcal D}}

\def\cf{{\mathcal F}}

\def\ch{{\mathcal H}}

\def\cn{{\mathcal N}}

\def\cs{{\mathcal S}}

\def\cz{{\mathcal Z}}

\def\ga{{\mathfrak A}} \def\gpa{{\mathfrak a}}
\def\gb{{\mathfrak B}}

\def\gf{{\mathfrak F}}
 \def\gpg{{\mathfrak g}}

 \def\gps{{\mathfrak s}}
 \def\gpt{{\mathfrak t}}

\def\gz{{\mathfrak Z}}

\def\bc{{\mathbb C}}

\def\bm{{\mathbb M}}
\def\bn{{\mathbb N}}

\def\br{{\mathbb R}}
\def\bt{{\mathbb T}}

\def\bz{{\mathbb Z}}

\def\a{\alpha}
\def\b{\beta}
\def\g{\gamma} \def\G{\Gamma}
\def\d{\delta}  
\def\eeps{\epsilon}
\def\eps{\varepsilon}
\def\l{\lambda} \def\L{\Lambda}

\def\m{\mu}

\def\n{\nu}
\def\r{\rho}
\def\s{\sigma} \def\S{\Sigma}
\def\t{\tau}
\def\f{\varphi}\def\ff{\phi} 
\def\th{\theta}  
\def\om{\omega} \def\Om{\Omega}
\def\z{\zeta}

\def\id{\hbox{id}}

\newtheorem{thm}{Theorem}[section]

\newtheorem{cor}[thm]{Corollary}

\newtheorem{rem}[thm]{Remark}

\def\aut{\mathop{\rm Aut}}
\def\ad{\mathop{\rm Ad}}

\def\sp{\mathop{\rm sp}}

\newcommand{\ty}[1]{\mathop{\rm {#1}}}
\def\di{\mathop{\rm d}\!}

\newcommand{\tr}{\mathop{\rm Tr}\!}
\newcommand{\nn}{\nonumber}

\begin{document}

\title[quantum disordered systems]
{some topics in quantum disordered systems}
\author{Stephen Dias Barreto}
\address{Stephen Dias Barreto\\
Department of Mathematics \\
Padre Conceicao College of Engineering\\
Verna Goa 403 722, India} \email{{\tt sbarreto@pccegoa.org}}
\author{Francesco Fidaleo}
\address{Francesco Fidaleo,
Dipartimento di Matematica, Universit\`{a} di Roma Tor Vergata,
Via della Ricerca Scientifica 1, Roma 00133, Italy} \email{{\tt
fidaleo@mat.uniroma2.it}}

\begin{abstract}
We discuss some recent results connected with the properties of
temperature states of quantum  disordered systems. This analysis
falls within the natural framework of operator algebras. Among the
results quoted here, we recall some ergodic and spectral
properties of KMS states, the possibility of defining the chemical
potential independently of the disorder, and finally the question
of the Gibbsianess for KMS states. This analysis can be considered
as a step towards fully understanding the very complicated
structure of the set of temperature states of quantum spin
glasses, and its connection with the breakdown of the symmetry for
replicas. \vskip 0.3cm \noindent
{\bf Mathematics Subject Classification}: 46L55, 82B44, 46L35.\\
{\bf Key words}: Non commutative dynamical systems; Disordered
systems; Classification of $C^{*}$--algebras, factors.
\end{abstract}

\maketitle

\section{introduction}
\label{sec1}

It is well--known that a disordered system, such as an alloy,
exhibits a very complicated thermodynamical behaviour. Recently,
the investigation of the structure of temperature states of spin
glasses, and the questions pertaining to the breakdown of symmetry
for replicas has generated considerable interest. Some interesting
connections with Archimedean number theory, and with
noncommutative geometry have been considered in order to
understand the above mentioned problems. Unfortunately, many
questions are still open. The reader is referred to \cite{B1, G, Ko, 
MPV, Ne, NeS1, PS} and the literature cited therein
for recent results along these lines, and for 
further details.

Another promising approach to study the complex behaviour of the
temperature states of a spin glass, is to use the standard
techniques of operator algebras, the last being very fruitful for
non disordered models, see e.g. \cite{BR1, BR2}. This approach
started with the seminal paper \cite{K}. There, besides some
general results concerning the properties of equilibrium states,
the phenomenon of weak Gibbsianess for quantum disordered systems
was pointed out. Recently, this strategy has been reconsidered,
see \cite{B, BF, F}.

In the present paper we discuss in an unified way such recent
results about disordered systems, including spin glasses. This
analysis covers states obtained by infinite volume limits of
finite volume Gibbs states, that is the so called quenched
disorder.

The paper is organized as follows. A preliminary section (Section
\ref{sec2}) contains the early definitions related to the
mathematical model of a quantum disordered system, and some
ergodic results connected with the well--known phenomenon of
self--averaging. Section \ref{sec3}, based on results of \cite{B,
BF}, is devoted to the structure of the set of the KMS states of
disordered systems, and some natural spectral properties. The
model describing the Edwards--Anderson quantum spin glass is
treated in some detail. Section \ref{sec4}, based on results in
\cite{F}, concerns the algebraic description for disordered
systems, of the chemical potential associated with the Gibbs
grand--canonical ensemble. Following the strategy of the pivotal
work \cite{AKTH}, it is shown that one can exhibit a
description of the chemical potential which does not depend on the
disorder. As for non disordered systems, it is shown that the
chemical potential is intrinsically connected with the
Connes--Radon--Nikodym cocycles (\cite{C, St}) associated with the
states obtained by composing the KMS state with the automorphisms
carrying the localizable (abelian) charges of the disordered
system under consideration. Section \ref{sec5} is based on the
paper \cite{K} and is devoted to the question of the joint Gibbsianess of a
KMS state of a disordered system.

The question of Gibbsianess arises when one considers infinite
volume limits of states obtained from the finite volume Gibbs
canonical ensemble on the skew space of joint configurations of
spins and couplings. The measure on such a skew space satisfies
the following properties. Its marginal distribution of the
couplings is the given probability measure describing the
disorder, whereas the conditional distribution of the spins, given
the couplings, is some infinite volume Gibbs state almost surely.
Such a field of infinite volume Gibbs states satisfies the
equivariance property \eqref{eqiv}, that is it gives an
Aizenman--Wehr metastate, see \cite{AW}. It can happen that the
measure so obtained is not necessarily jointly Gibbsian. This
issue was addressed in quantum setting in Proposition 5.14 of
\cite{K}, and it is checked for classical disordered models in
\cite{EMSS, Ku}. In Section \ref{sec5} we present a shortened
discussion of these questions.

To end the introduction, we point out that a large number of
problems remain open in the study of disordered systems. However,
the investigations of disordered systems by the standard
techniques of operator algebras could be a significant step
towards fully understanding the very complicated structure of the
set of the KMS states of quantum spin glasses, and its connection
with the breakdown of the symmetry for replicas.

\section{preliminaries}
\label{sec2}

Our starting point for the investigation of the disordered systems
is a separable $C^{*}$--algebra $\ca$ with an identity $I$,
describing the physical observables. Sometimes (i.e. Section
\ref{sec4}), $\ca$ can be obtained as the fixed point algebra
$\ca=\cf^{G}$ under a strongly continuous action
$$
\g:g\in G\mapsto\g_{g}\in\aut(\cf)
$$
of a compact second countable group $G$ (the {\it gauge group}) on
another separable $C^{*}$--algebra $\cf$ (the {\it field
algebra}), see e.g. \cite{DHR1, DHR2, DR1, DR2}. We suppose that
the group $\{\a_{x}\}_{x\in\bz^{d}}$ of spatial translations acts
on $\cf$.\footnote{To simplify the matter, we are assuming that
$\bz^{d}$ is the group of the space symmetries. Some of the main
results quoted below can be generalized to ``continuous''
disordered systems, where the group of the space symmetries is
$\br^{d}$.} The sample space for the couplings is a standard
measure space $(X,\n)$, $X$ being a compact separable space, and
$\n$ a Borel probability measure. The group $\bz^{d}$ of the
spatial translations is supposed to act on the probability space
$(X,\n)$ by measure preserving transformations
$\{T_{x}\}_{x\in\bz^{d}}$.

A one parameter random group of automorphisms
\begin{equation*}
(t,\xi)\in\br\times X\mapsto\t^{\xi}_{t}\in\aut(\cf)
\end{equation*}
is acting on $\cf$. It is supposed to be strongly continuous in
the time variable for each fixed $\xi\in X$, and jointly strongly
measurable. We assume that $\t$ acts locally, that is the strongly
measurable, essentially bounded function
$f_{A,t}(\xi):=\t^{\xi}_{t}(A)$ belongs to the $C^{*}$--subalgebra
$\cf\otimes L^{\infty}(X,\n)$, where the above $C^{*}$--tensor
product is uniquely determined as any commutative $C^{*}$--algebra
is nuclear. Furthermore, we assume that the random time evolution
described above is nontrivial, if it is not otherwise
specified.\footnote{This nontriviality condition simply means that
the one parameter group $\{\gpt_{t}\}_{t\in\br}$ given in
\eqref{aztemp1} is nontrivial.} We assume the commutation rule
\begin{equation}
\label{4} \t^{T_{x}\xi}_{t}\a_{x}=\a_{x}\t^{\xi}_{t}
\end{equation}
for each $x\in\bz^{d}$, $\xi\in X$, $t\in\br$. If we start
directly with the observable algebra, we assume analogous
properties described above for $\ca$ instead of $\cf$.

If we start with the field algebra $\cf$, we suppose further that
\begin{equation*}
\a_{x}\g_{g}=\g_{g}\a_{x}\,,
\t^{\xi}_{t}\g_{g}=\g_{g}\t^{\xi}_{t}\,,
\end{equation*}
for each $x\in\bz^{d}$, $\xi\in X$, $t\in\br$, and $g\in G$. If
$\ca$ is obtained by a principle of gauge invariance, then, by the
above commutation rules, $\a_{x}$ and $\t^{\xi}_{t}$ leave $\ca$
globally stable.

We address also the situation when Fermion operators are present
in $\cf$. Namely, there exists an automorphism $\s$ of $\cf$
commuting with the all the gauge transformations, the spatial
translations and the random time evolution, such that $\s^{2}=e$.
We put
\begin{equation}
\label{fesi} 
\cf_{+}:=\frac{1}{2}(e+\s)(\cf)\,,\quad
\cf_{-}:=\frac{1}{2}(e-\s)(\cf)\,.
\end{equation}

In order to achieve the disorder, it is natural to set for the
field algebra,
$$
\gf:=\cf\otimes L^{\infty}(X,\n)\,.
$$

Notice that, by identifying $\gf$ with a closed subspace of
$L^{\infty}(X,\n;\cf)$, each element $A\in\gf$ is uniquely
represented by a measurable essentially bounded function
$\xi\mapsto A(\xi)$ with values in $\cf$.
We have on $\gf$,
\begin{align}
\label{aztemp1}
&\gpa_{x}(A)(\xi):=\a_{x}(A(T_{-x}\xi))\,,\nn\\
&\gpt_{t}(A)(\xi):=\t_{t}^{\xi}(A(\xi))\,,\nn\\
&\gpg:=\g\otimes\id_{L^{\infty}(X,\n)}\,,\\
&\gps:=\s\otimes\id_{L^{\infty}(X,\n)}\,.\nn
\end{align}

It is easy to verify that $\{\gpa_{x}\}_{x\in\bz^{d}}$,
$\{\gpt_{t}\}_{t\in\br}$ and $\{\gpg_{g}\}_{g\in G}$ define
actions of $\bz^{d}$, $\br$ and $G$ on $\gf$ which are mutually
commuting, and commute also with the parity automorphism $\gps$.
The subspaces $\gf_{+}$ and $\gf_{-}$ are defined as in
\eqref{fesi}. Furthermore, $\gpa_{x}$ and $\gpt_{t}$ leave the
disordered observable algebra $\ga$ globally stable. Namely,
$\{\gpa_{x}\}_{x\in\bz^{d}}$ and $\{\gpt_{t}\}_{t\in\br}$ define
by restriction, mutually commuting actions of $\bz^{d}$ and $\br$
on $\ga$, respectively.

In order to study a class of states of interest for disordered
systems, we start with  $*$--weak measurable fields of states
$$
\xi\in X\mapsto\f_{\xi}\in\cs(\cf)\,.
$$

We suppose that the field $\{\f_{\xi}\}_{\xi\in X}$ fulfills
almost surely, the equivariance condition
\begin{equation}
\label{eqiv} \f_{\xi}\circ\a_{x}=\f_{T_{-x}\xi}
\end{equation}
w.r.t. the spatial translations, simultaneously.
A state $\f$ on $\gf$ is naturally defined as
\begin{equation}
\label{ddf} \f(A)=\int_{X}\f_{\xi}(A(\xi))\n(\di\xi)\,, \qquad
A\in\gf\,.
\end{equation}

It is immediate to verify that $\f$ defined as above is invariant
w.r.t. the space translations $\gpa_{x}$. Moreover,
$\f\lceil_{I\otimes L^{\infty}(X,\n)}$ is a normal state.

Equally well, one can start with a  $\gpa$--invariant state $\f$
on $\gf$, which is normal when restricted to $I\otimes
L^{\infty}(X,\n)$ (indeed, $\f\lceil_{I\otimes
L^{\infty}(X,\n)}=\int_{X}\,\cdot\,\n(\di\xi)$). Then, we can
recover a $*$--weak measurable field $\{\f_{\xi}\}_{\xi\in
X}\subset\cf$ fulfilling \eqref{eqiv}. Such a measurable field
provides the direct integral decomposition of $\f$ as in
\eqref{ddf}, see \cite{BF}, Theorem 4.1, see also \cite{K}, 
Proposition 4.1. Similar considerations
can be applied to the observable algebras $\ga$ as well. In the
sequel, we denote by $\cs_{0}(\ga)$, $\cs_{0}(\gf)$ the convex
closed subsets of states on $\ga$, $\gf$ respectively, fulfilling
the properties listed above.

Now, we recall same properties of states on $\gf$ or
$\ga$. We restrict ourselves to the field algebra, the other case
being similar. Let $C,D\in\gf$, and $A,B\in\gf_{+}\bigcup\gf_{-}$.
Put $\eeps_{A,B}=-1$ if $A,B\in\gf_{-}$ and $\eeps_{A,B}=1$ in the
three remaining possibilities. We say that the state
$\f\in\cs(\gf)$ is {\it asymptotically Abelian} w.r.t. $\gpa$ if
\begin{equation*}
\lim_{|x|\to+\infty}
\f\left(C\big(\gpa_{x}(A)B-\eeps_{A,B}B\gpa_{x}(A)\big)D\right)=0\,,
\end{equation*}

The state $\f\in\cs(\gf)$ is {\it weakly clustering} w.r.t. $\gpa$
if
\begin{equation*}
\lim_{N}\frac{1}{|\L_{N}|}\sum_{x\in\L_{N}}
\f(A\gpa_{x}(B))=\f(A)\f(B)\,,
\end{equation*}
$\L_{N}$ being the box with a vertex at the origin, containing
$N^{d}$ points with positive coordinates. Notice that, Many of
interesting states arising from quantum physics are naturally
asymptotically Abelian w.r.t. the spatial translations, see e.g.
\cite{L}, see also \cite{BF}, Proposition 2.3 for the case of
disordered systems.

It is well--know that for a $\gpa$--invariant asymptotically
Abelian state, the $\gpa$--weak clustering property is equivalent
to the $\gpa$--ergodicity, see e.g. Proposition 5.4.23 of
\cite{BR2}. We speak, without any further mention, and if it is
not otherwise specified, about asymptotic Abelianess, weak
clustering, or ergodicity for states, if they satisfy these
properties w.r.t. the spatial translations.

In \cite{F}, Section 3, it is shown under the ergodicity of
the action $\{T_{x}\}_{x\in\bz^{d}}$ on the sample space $(X,\n)$, 
that $\f\in\cs_{0}(\gf)$ is weakly
clustering if and only if
\begin{equation*}
\lim_{N}\frac{1}{|\L_{N}|}\sum_{x\in\L_{N}}
\f_{\xi}(A\a_{x}(B(T_{-x}\xi)))=\f_{\xi}(A)\f(B)\,,\quad
A\in\cf\,,B\in\gf
\end{equation*}
in the weak  topology of $L^{1}(X,\n)$, and if
\begin{equation*}
\lim_{|x|\to+\infty}\f_{\xi}(A\a_{x}(B))=\f_{\xi}(A)\f_{T_{-x}\xi}(B)\,,\quad
A,B\in\cf
\end{equation*}
almost surely. Here, $\{\f_{\xi}\}_{\xi\in X}\subset\cf$ is the
measurable field of states associated with $\f$, satisfying
\eqref{eqiv}. The last results are connected with the well--known property
of {\it self--averaging} for disordered systems.

For the definition of the KMS boundary condition, the equivalent
characterizations of the KMS boundary condition, the main results
about KMS states, and finally the connections with Tomita theory
of von Neumann algebras, we refer the reader to \cite{BR2, St} and
the references cited therein. We have for the modular group
$\s^{\ff}$ of a KMS state $\ff$,
\begin{equation}
\label{modgns} \s^{\ff}_{t}\circ\pi_{\ff}=\pi_{\ff}\circ\t_{-\b
t}\,,
\end{equation}
$\pi_{\ff}$ being the GNS representation of $\ff$.
\begin{rem}
\label{rm11} It is a well--known fact that $\f\in\cs_{0}(\gf)$ is
a KMS state if and only if $\{\f_{\xi}\}_{\xi\in X}$ are KMS
states almost surely, see e.g. \cite{F, S}.
\end{rem}

\section{the structure of the KMS states and their spectral properties}
\label{sec3}

The present section is based on results contained in \cite{B, BF},
to which the reader is referred for further details. To simplify the 
matter, we suppose that the algebra of observables $\ca$ is a
quasi--local algebra (see e.g. \cite{BR1}, Section 2.6) whose
local algebras are isomorphic to full matrix algebras, even if
most of the forthcoming analysis works for more general
situations. We further suppose that the action $x\mapsto T_{x}$ of
the spatial translations on the sample space $(X,\n)$ is ergodic.

The first result concerns the almost surely independence of the
Arveson spectrum for random systems.\footnote{The Arveson
spectrum, known in physics as the set of the {\it Bohr
frequencies} is made up of frequencies of the quanta which a
physical system (like an atom) may emit or absorb, see \cite{FM},
Theorem 5.3 for a simple example. For the definition of the Arveson
spectrum, see e.g. \cite{P}.} This is precisely Theorem 5.3 of
\cite{B} (see also \cite{BF}, Theorem 2.1). For the reader's
convenience, we report here its proof.
\begin{thm}
\label{ba} Under the above assumptions, there exists a measurable
set $F\subset X$ of full measure, and a closed set $\S\subset\br$
such that $\xi\in F$ implies $\sp(\t^{\xi})=\S$.
\end{thm}
\begin{proof}
We get by \cite{P}, Proposition 8.1.9,
$$
\sp(\t^{\xi})=\bigcap_{f\in L^{1}(\br)}
\big\{s\in\br\;\big|\;|\hat{f}(s)|\leq\|\t^{\xi}_{f}\|\big\}
$$
where $\hat{f}$ is the inverse Fourier transform of $f$, and
$$
\t^{\xi}_{f}(A):=\int_{-\infty}^{+\infty}f(t)\t^{\xi}_{t}(A)\di
t\,.
$$

By a standard density argument, we can reduce the situation to a
dense set $\{f_{k}\}_{k\in\bn}\subset L^{1}(\br)$. Define
$\G_{k}(\xi):=\|\t^{\xi}_{f_{k}}\|$. It was shown in \cite{B} that
the functions $\G_{k}$ are measurable and invariant. By
ergodicity, they are constant almost everywhere. Let
$\{N_{k}\}_{k\in\bn}$ be null subsets of $X$ such that, for each
$k\in\bn$ and $\xi\in N_{k}^{{}^{c}}$,
$$
\G_{k}(\xi)=\|\G_{k}\|_{\infty}\,.
$$

Consider $F:=\big(\bigcup_{k\in\bn}N_{k}\big)^{c}$, and take
$\S:=\sp(\t^{\xi_{0}})$, where $\xi_{0}$ is any element of $F$. We
have that $F$ is a measurable set of full measure, and $\xi\in F$
implies $\sp(\t^{\xi})=\S$.
\end{proof}

Let $\om\in\cs_{0}(\ga)$, $\{\om_{\xi}\}_{\xi\in X}\subset\ca$ be
the measurable field of states associated to $\om$ as described
above, and $\pi_{\om}$, $\pi_{\om_{\xi}}$ the GNS representations
of $\om$, $\om_{\xi}$ respectively. Consider the subcentral
decomposition
\begin{equation}
\label{11} M=\int^{\oplus}_{X}M_{\xi}\,\n(\di\xi)
\end{equation}
of $M:=\pi_{\om}(\ga)''$. It is shown in Section 3 of \cite{BF},
that $M_{\xi}$ coincides with $\pi_{\om_{\xi}}(\ca)''$ almost
surely.

Some of the main results in \cite{BF} (see Section 4 of that
paper) are collected in the following theorem. Here, we sketch its
proof for the convenience of the reader.
\begin{thm}
\label{main} 
Let $\om\in\cs_{0}(\ga)$ be a KMS state at inverse
temperature $\b\neq0$. Then only type $\ty{III_{\l}}$ factors,
$\l\in(0,1]$, can appear in the central decomposition of
$\pi_{\om}(\ga)''$.

In addition, if $\gz_{\pi_{\om}}\sim L^{\infty}(X,\n)$, then there
exists a unique $\l\in(0,1]$ such that $\pi_{\om_{\xi}}(\ca)''$
are type $\ty{III_{\l}}$ factors almost surely.
\end{thm}
\begin{proof}
Taking into account \eqref{modgns} and Theorem \ref{ba}, it is
shown in Proposition 4.5 of \cite{BF} that
$$
\G_{B}(\pi_{\f}(\ga)'')\equiv-\b\sp(\gpt)=-\b\sp(\t^{\om})
$$
almost surely, where $\G_{B}$ denotes the Borchers invariant
(\cite{Bo}) and ``$\sp$'' the Arveson spectrum.\footnote{This is a
consequence of general results of \cite{A0, HL}.} This means that,
in our situation, the infinite semifinite portion, and
$\ty{III_{0}}$ portion in the factor decomposition of
$\pi_{\om_{\xi}}(\ca)''$ are avoided (cf. Theorem 4.4 of
\cite{BF}). The $\ty{I_{\text{fin}}}$ is trivially avoided, and
finally the $\ty{II_{1}}$ is avoided as we are assuming that
$\gpt_{t}$ is nontrivial (i.e. $\om$ is not a trace). This proves
the first part.

Now, $\gz_{\pi_{\om}}\sim L^{\infty}(X,\n)$ means that the direct
integral decomposition \eqref{11} is the factor decomposition of
$\pi_{\om}(\ga)''$. In this situation,
$\G_{B}(\pi_{\om}(\ga)'')=\G(\pi_{\om_{\xi}}(\ca)'')$ almost
surely, where $\G$ is Connes $\G$--invariant (\cite{C}). This
implies by the previous assertion, that there exists a
$\l\in(0,1]$ such that $\pi_{\om_{\xi}}(\ca)''$ is a type
$\ty{III_{\l}}$ factor almost surely.
\end{proof}
\begin{rem}
It is explained in \cite{BF} that states $\om\in\cs_{0}(\ga)$ for
which $\gz_{\pi_{\om}}\sim L^{\infty}(X,\n)$ can be thought as
describing the ``pure termodynamical phase'' in the case of
disordered systems.\footnote{For classical disordered systems, the
centre $\gz_{\pi_{\om}}$ should be replaced by {\it the algebra at
infinity}
$$
\gz^{\perp}_{\om}:=\bigwedge_{\L\,\text{finite}}\pi_{\om}(\ga_{\L^{c}})''\,.
$$}
The fact that $\pi_{\om_{\xi}}$ generates a type $\ty{III_{\l}}$
von Neumann algebra, $\l\in(0,1]$, is in accordance with the
standard fact that the physically relevant quantities do not
depend on the disorder.
\end{rem}

In order to demonstrate some possible applications, we specialize
our analysis to the quantum model described by the formal random
Hamiltonian
\begin{equation}
\label{fham}
H=-\frac{1}{2}\sum_{i\in\bz^{d}}\sum_{|i-j|=1}J_{i,j}\s(i)\s(j)\,.
\end{equation}

Here, $\s(i)$ is the spin--operator along the $z$--axis (i.e. the
Pauli matrix $\s_{z}$) acting on the $i$--th site, and the
$J_{i,j}$ are independent identically distributed random variables
with common distribution $p(\di y)$ on the real line. To simplify,
we suppose that the law $p$ is compactly supported.\footnote{The
model with symmetric $p$ is known in literature as the (quantum)
{\it Edwards--Anderson spin glass}, see \cite{EA}.}

We first suppose that for a fixed $\b>0$, the Ising--type model
described by the Hamiltonian \eqref{fham} admits a unique KMS
state, say $\om_{\xi}$, almost surely. Then the map $\xi\in
X\mapsto \om_{\xi}\in\cs(\ca)$ is automatically $*$--weak
measurable and satisfies almost surely the
equivariance condition  \eqref{eqiv}.\footnote{This situation
arises if the quantum model under consideration admits some
critical temperature. The situation is well clarified for many
classical disordered models (see e.g. \cite{Ne}), contrary to the
quantum situation where, in the knowledge of the authors, there
are very few rigorous results concerning this fact. However, it is
expected that quantum disordered systems also exhibit critical
temperatures in the high temperature regime.} In this situation,
there exists a unique state $\om\in\cs_{0}(\ga)$ for the model
under consideration, given by
\begin{equation}
\label{cio} 
\om(A)=\int_{X}\om_{\xi}(A(\xi))\n(\di\xi)\,,\qquad
A\in\ga\,.
\end{equation}

Furthermore, each of the other states $\om_{f}\in\cs(\ga)$
which are normal when restricted to $I\otimes L^{\infty}(X,\n)$,
are in one--to--one correspondence with the positive normalized
functions  $f\in L^{\infty}(X,\n)$ through the relation
\begin{equation}
\label{cio1}
\om_{f}(A)=\int_{X}f(\xi)\om_{\xi}(A(\xi))\n(\di\xi)\,,\qquad
A\in\ga\,.
\end{equation}
Taking into account Theorem \ref{main}, it follows that KMS states
in \eqref{cio}, \eqref{cio1} generate type $\ty{III_{\l}}$ von
Neumann algebras for a fixed $\l\in(0,1]$ depending on the
temperature.

We end the present section by briefly describing the multiple
phase regime. After taking the infinite--volume limit along
various subsequences $\L_{n_{k}}\uparrow\bz^{d}$, we will find by
compactness, in general different $\gpt$--KMS states in
$\cs_{0}(\ga)$ at fixed inverse temperature $\b$. Fix one such a
state $\psi\in\cs_{0}(\ga)$, with the associated $*$--weak
measurable field $\xi\in X\mapsto\psi_{\xi}\in\cs(\ca)$ of
$\t^{\xi}$--KMS states respectively, satisfying the equivariance
property \eqref{eqiv}. According to Proposition 3.1 of \cite{BF},
the set of the $\gpt$--KMS states $\psi_{T}\in\cs(\ga)$, normal
w.r.t. $\psi$, have the form
$$
\psi_{T}(A)=\int_{X}
\big\langle\pi_{\psi_{\xi}}(A(\xi))T(\xi)^{1/2}\Om_{\psi_{\xi}},
T(\xi)^{1/2}\Om_{\psi_{\xi}}\big\rangle_{\ch_{\psi_{\xi}}}\n(\di\xi)\,.
$$

Here, $(\pi_{\psi_{\xi}},\ch_{\psi_{\xi}},\Om_{\psi_{\xi}})$ is
the GNS representation of $\psi_{\xi}$, $\{T(\xi)\}_{\xi\in X}$ is
a measurable field of closed densely defined operators on
$\ch_{\psi_{\xi}}$ affiliated to the (isomorphic) centres
$\gz_{\psi_{\xi}}$ respectively, satisfying
$\Om_{\psi_{\xi}}\in\cd_{T(\xi)^{1/2}}$ almost surely, and
${\displaystyle
\int_{X}\|T(\xi)^{1/2}\Om_{\psi_{\xi}}\|_{\ch_{\psi_{\xi}}}^{2}\n(\di\xi)=1}$.

\section{an algebraic description of the chemical potential} 
\label{sec4}

The present section is based on results contained in \cite{F}, to
which the reader is referred for further details. As in the
previous section, we suppose that the action $x\mapsto T_{x}$ of
the spatial translations on the sample space $(X,\n)$ is ergodic.
Here, in order to give an algebraic description of the chemical
potential, we suppose that the algebra of the observables $\ca$ is
obtained by a principle of gauge invariance, as the fixed point
subalgebra of the field algebra $\cf$.

The technical problem which is behind the algebraic description of
the chemical potential is the question of extending a KMS state on
the observable algebra to a KMS state on the field algebra. This
is a very delicate issue studied in a series of papers connected
with non disordered models, see \cite{AKTH, AK, KT1, KT2, L1}. The
strategy and the main results of the pivotal paper \cite{AKTH} can
be effectively used to investigate this problem for disordered
systems.

After proving some preliminary but crucial results such as the
facts that the stabilizer and the asymmetry subgroup of a state
$\f\in\cs_{0}(\gf)$ coincide almost surely with the corresponding
objects associated with  $\f_{\xi}$,\footnote{The stabilizer
$G_{\f}$ of $\f\in\cs(\gf)$ is defined as ${\displaystyle G_{\f}:=
\big\{g\in G\,\big|\,\f\circ\gpg_{g}=\f\big\}}$. For
$\f\in\cs_{0}(\gf)$, the stabilizers $G_{\f_{\xi}}$ of the
$\f_{\xi}\in\cs(\cf)$ associated with $\f$, are defined
analogously. We denote by an abuse of notation, the normalizer and
the centralizer of a subgroup $H\subset G$ by $\cn(H)$ and
$\cz(H)$ respectively, and are defined as ${\displaystyle\cn(H):=
\big\{g\in G\,\big|\,gHg^{-1}=H\big\}}$,
${\displaystyle\cz(H):=\big\{g\in G\,\big|\,gh=hg\,,h\in
H\big\}}$. The asymmetry subgroup of an invariant state is the
subgroup for which the spectrum of the associated unitary
implementation is one--sided, see \cite{AKTH}, Definition II.3.}
it is shown in \cite{F} that the KMS states in $\cs_{0}(\ga)$ satisfying certain
natural properties, can be extended to KMS states on
$\cs_{0}(\gf)$.

We report here the main result of \cite{F} concerning the
algebraic description of the chemical potential for disordered
systems.
\begin{thm}
\label{gaga2} Let $\f\in\cs_{0}(\gf)$ be a weakly clustering
asymptotically Abelian state whose restriction to $\ga$ is
$(\gpt,\b)$--KMS state at inverse temperature $\b\neq0$.

Then there exist a closed subgroup $N\subset G_{\f}$, a continuous
one parameter subgroup $t\in\br\mapsto\eps_{t}\in\cz(G_{\f})$, a
continuous one parameter subgroup $t\in\br\mapsto\z_{t}\in
G_{\f}$, and a measurable subset $F\subset X$ of full measure such
that, for each $\xi\in F$,
\begin{itemize}
\item[(i)] the $N$--spectrum of $\f_{\xi}$ is one--sided,
\item[(ii)] the restriction of $\f_{\xi}$ to
$\cf^{N}:=\big\{A\in\cf\,\big|\,\g_{g}(A)=A\,,g\in N\big\}$ is a
$(\th^{\xi},\b)$--KMS state for the modified time evolution
$\th^{\xi}_{t}:=\t^{\xi}_{t}\g_{\eps_{t}\z_{t}}$,
\item[(iii)] the image $[\z_{t}]:=\z_{t}N$ in $G_{\f}/N$ is in
$\cz(G_{\f}/N)$.
\end{itemize}
\end{thm}
The proof is a consequence of Theorem II.4 of \cite{AKTH}, taking
into account the above mentioned results pertaining to the
stablizers and the asymmetry subgroups, see \cite{F}, Theorem 4.7
for further details.

Theorem \ref{gaga2} can be directly applied to the case $\b=0$
which corresponds to the case when the restriction to $\ga$ of a
state in $\cs_{0}(\gf)$ is a trace.
\begin{cor}
Under the same assumption of Theorem \ref{gaga2} but for $\b=0$,
there exist a closed subgroup $N\subset G_{\f}$, a continuous one
parameter subgroup $t\in\br\mapsto\z_{t}\in G_{\f}$, and a
measurable subset $F\subset X$ of full measure such that the
assertions of Theorem \ref{gaga2} hold true with $\rm{(ii)}$
replaced by
\begin{itemize}
\item[(ii')] the restriction of $\f_{\xi}$ to
$\cf^{N}:=\big\{A\in\cf\,\big|\,\g_{g}(A)=A\,,g\in N\big\}$ is a
$(\g_{\z_{t}},-1)$--KMS state.\footnote{The value $-1$ is chosen
such that the restriction of the modular group to $\pi_{\f}(\gf)$
coincides with the time evolution $\g_{\z_{t}}$, see \eqref{modgns}.}
\end{itemize}
\end{cor}
\begin{proof}
We start by noticing that a state $\om$ satisfies the KMS boundary
condition at inverse temperature $0$ if and only if it is a
trace.\footnote{If $\om\in\cs_{0}(\ga)$ is a trace, then
$\om_{\xi}$ is a trace almost surely, the proof being the same as
that of Proposition 3.2 of \cite{BF}.} This means that $\om$ is
KMS at any inverse temperature $\b\neq0$ w.r.t. the trivial one
parameter automorphism group $t\mapsto\id_{\ga}$. By convention,
we choose $\b=-1$.

The proof follows by the previous theorem, taking into account
that, in this situation, we can choose for
$t\in\br\mapsto\eps_{t}\in\cz(G_{\f})$ the trivial one parameter
automorphism group, see \cite{AKTH}, Theorem 2.2.
\end{proof}

Now, we apply the previous results to a simple case in order to
understand how the chemical potential appears. To avoid technical
complications, we further assume that $\cf$ is a quasi--local
algebra whose local algebras are isomorphic to full matrix
algebras. We assume also that the gauge group is the unit circle
$\bt$.\footnote{One can directly start from the observable algebra
$\ca$, and reconstruct the field algebra $\cf$ from the
localizable charges of interest of the model, see \cite{DHR2,
DR2}. See also \cite{NS} for a possible application of this
strategy to the situation considered here.} In this situation, the
localizable charges of the model under consideration are
associated with the powers $[\s^{n}]\in\text{Out}(\ca)$ of a
single localized transportable outer automorphism $\s$, see
\cite{DHR1}.

Let $\om\in\cs_{0}(\ga)$ be a $(\gpt,\b)$--KMS state at inverse
temperature $\b\neq0$ such that $\gz_{\pi_{\om}}\sim
L^{\infty}(X,\n)$. Take a weakly clustering extension $\f$ of
$\om$ to all of $\gf$ which exists by Proposition 4.2 of
\cite{BF}. In order to avoid the possibility of null chemical
potential, we further suppose that $\f$ is gauge invariant. In this
situation, the asymmetry subgroup $N_{\f}$ of $\f$ is trivial, or
equivalently $N_{\f_{\xi}}$ is trivial 
almost surely.\footnote{In this situation, we have 
$\sp(U_{\f})=\bz\equiv\widehat{G}$. Indeed, let $V_{\s}\in\cf$ be the unitary 
implementing the automorphism $\s$ on $\ca$. This means that 
$\g_{\th}(V_{\s})=e^{i\th}V_{\s}$. The vector 
$\Psi_{n}:=\pi_{\f}(V_{\s}^{n}\otimes I)\Om_{\f}$ is an eigenvector 
of $U_{\f}(\th)$ corresponding to the eigenvalue $e^{in\th}$.}

Let $\r$ be a localized automorphism of $\ca$ carrying the charge
$n$ (i.e. $\r\in[\s^{n}]$), and $\om\in\cs_{0}(\ga)$ as above.
Then, $\r$ extends to a measurable field $\{\r_{\xi}\}_{\xi\in X}$
of normal automorphisms of the weak closure
$\pi_{\om_{\xi}}(\ca)''$, almost surely, see \cite{F}, Proposition
5.1.

Consider the unitary $U$ implementing $\r$ on $\ca$, together with
the state $\f_{U}:=\f\circ\ad(U\otimes I)$. We have for the
Connes--Radon--Nikodym cocycle (\cite{C, St}),
\begin{align}
\label{crncoc} \big(D\f_{U}:D\f)_{t}=&\int^{\oplus}_{X}
\pi_{\f_{\xi}}(U^{*})\s^{\f_{\xi}}_{t}(\pi_{\f_{\xi}}(U))\n(\di\xi)\nn\\
=e^{in\b\m t}&\int^{\oplus}_{X} \pi_{\f_{\xi}}(U^{*}\t^{\xi}_{-\b
t}(U))\n(\di\xi)\,,
\end{align}
for some $\m\in\br$.\footnote{Notice that the direct integral decomposition 
of $\big(D\f_{U}:D\f)$ exists by
the results in Appendix A of \cite{BF}, taking into account the
definition of the associated modular operator in terms of the
balanced weight, see e.g. \cite{St}, Section 3.} 

Here, we have used
$\g_{\th}(U)=e^{in\th}U$, and
$\s^{\f}_{t}\circ\pi_{\f}=\pi_{\f}\circ\gpt_{-\b t}\circ\gpg_{\b\m
t}$ by Theorem \ref{gaga2}, taking into account \eqref{modgns}.

Now, we take advantage of the fact that $\om\circ(\r\otimes\id)$
extends to a normal state on all of $\pi_{\om}(\ga)''$. Denote by
the same symbol $\om_{\xi}$ the normal extension of $\om_{\xi}$
itself to all of $\pi_{\om_{\xi}}(\ca)''$. Under this
identification, the equivariant measurable field of states
$\{\om_{\xi}\circ\r_{\xi}\}_{\xi\in X}$ provides the direct
integral decomposition of the mentioned normal extension of
$\om\circ(\r\otimes\id)$. Here, the $\r_{\xi}$ are the above
mentioned normal automorphisms of $\pi_{\f_{\xi}}(\ca)''$
extending $\r$. We have by \eqref{crncoc} and the fact that
$U^{*}\t^{\xi}_{-\b t}(U)$ is gauge invariant,
\begin{equation}
\label{crncoc1} \big(D(\om_{\xi}\circ\r_{\xi}):D\om_{\xi})_{t}=
e^{in\b\m t} \pi_{\om_{\xi}}(U^{*}\t^{\xi}_{-\b t}(U))\
\end{equation}
almost everywhere.

Formula \eqref{crncoc1} explains the occurrence of the chemical
potential $\m\in\br$ as an object intrinically associated to the
observable algebra. Furthermore, according to this description, it
does not depend on the disorder for KMS states $\om\in\cs_{0}(\ga)$ such
that $\gz_{\pi_{\om}}\sim L^{\infty}(X,\n)$. This is in accordance
with the standard fact that the physically relevant quantities
should not depend on the disorder.

We end by noticing that the last analysis can be straightforwardly
extended to the $n$--dimensional torus or, more generally, to non
Abelian $n$--dimensional Lie groups, obtaining a $n$--parametric 
chemical potential. Furthermore, the results of
the present section could be applied to continuous disordered
systems such as possible disordered systems arising from quantum
field theory, even if some technical gaps need to be covered in
the last situation.

\section{the question of Gibbsianess}
\label{sec5}

The question of Gibbsianess arises when one takes infinite volume
limits of finite volume Gibbs states. If one considers such an
object as a state on the skew space of joint configurations of
spins and couplings, it can happen that the state so obtained is
not jointly Gibbsian.\footnote{Such a weakly Gibbsian state is
known in literature as a {\it metastate}, see \cite{AW} for the
classical situation.}

The possible appearance of weakly Gibbsian non jointly Gibbsian
states was firstly pointed out by Kishimoto in relation to the
quantum case, see \cite{K}, see also \cite{A1}. Recently, simple
examples of weakly Gibbsian states were constructed for classical
disordered systems, see \cite{EMSS, Ku}. In the present section we
discuss some results of \cite{K} connected with Gibbsianess.

We start by setting
$$
\ca:=\bigotimes_{\bz^{d}}\bm_{n}(\bc)\,.
$$

We assume that also the sample space has a local structure, and is
finite for finite regions. More precisely,
$$
X:=\prod_{\bz^{d}}{\bf k}\,.
$$

In this simple situation, we can take for $\ga$, the separable
$C^{*}$--algebra consisting of all the continuous functions from
$X$ with values in $\ca$:
$$
\ga:=C(X,\ca)\,.
$$

The algebra $\ga$ of observables is equipped with a natural local structure
${\displaystyle\big\{\ga_{\L}\,\big|\,\L\subset\bz^{d}\,,\L\, 
\text{bounded}\big\}}$,
with $\ga_{\L}=C(X_{\L},\ca_{\L})$.

It is immediate to show that we can recover $\ga$ from a principle
of gauge invariance from
$$
\gb:=
\bigotimes_{\bz^{d}}\big(\bm_{n}(\bc)\otimes\bm_{k}(\bc)\big)\,,
$$
equipped with the local structure
${\displaystyle\big\{\gb_{\L}\,\big|\,\L\subset\bz^{d}\,,\L\, 
\text{bounded}\big\}}$,
where 
$$
\ga_{\L}=
\bigotimes_{\L}\big(\bm_{n}(\bc)\otimes\bm_{k}(\bc)\big)\,.
$$

The gauge group is given by
$$
G:=\prod_{\bz^{d}}\bz_{k}\,,
$$
$\bz_{k}=\{0,1,\dots,k-1\}$ being the cyclic group of order $k$. Here, the gauge
action $a_{g}$, $g=\{l_{r}\}_{r\in\bz^{d}}$, is given by
\begin{equation}
\label{gaac} 
a_{\{l_{r}\}}:=
\id_{\ca}\bigotimes\bigg(\otimes_{r\in\bz^{d}}\ad(V^{l_{r}})\bigg)\,,
\end{equation}
with
$$
V:=\sum_{j=0}^{k-1}e^{i2\pi \frac{j}{k}}P_{j}\,,
$$
$\{P_{j}\}_{j=0,\dots,k-1}$ being the set of the minimal
projection of $\bc^{k}\subset\bm_{k}(\bc)$.\footnote{The
algebra $\ga$ slightly differs from the analogous
object in the previous sections without affecting the forthcoming
analysis. However, the principle of gauge invariance described in
the present section differs significantly from that assumed in
Section \ref{sec4}. Namely, the algebra $\gb$ is non disordered
(i.e. it is simple, contrary to $\gf$). In addition, $a_{g}$ in
\eqref{gaac} acts in a different way from $\gpg_{g}$ in
\eqref{aztemp1}. Finally, the gauge action \eqref{gaac} does not
commute with the space translations. Thus, the theory of the
chemical potential described in Section \ref{sec4} cannot be
directly applied to the field system $(\gb,G, a_{g})$.} This
situation covers the example described by the Hamiltonian
\eqref{fham} if the common distribution $p(\di y)$ for the
couplings has finite support made of $k$ elements, and falls within
the framework of \cite{K}, even if slightly more general
situations are treated there.

Let $\{\gpt_{t}\}_{t\in\br}$ be a strongly continuous one
parameter group of automorphisms of $\gb$ describing the time
evolution commuting with the space translations. We suppose that
$$
\gpt_{t}(\ga)=\ga\,,\quad t\in\br\,,
$$
that is $\gpt_{t}$ leaves globally stable $\ga$, and that
$$
\bigcup_{\L\,{\rm finite}}\gb_{\L}\subset\cd(\d)\,,
$$
$\cd(\d)$ being the domain of the infinitesimal generator of
$\gpt_{t}$.

Consider, for $\xi\in X$ the homomorphism
$$
\chi_{\xi}(A):=A(\xi)\,,\quad A\in\ga
$$
of $\ga$ onto $\ca$. It follows by Proposition 3.1 of \cite{K}
that a random time evolution $\{\t^{\xi}_{t}\}_{t\in\br}$ is
recovered as
$$
\t^{\xi}_{t}(\chi_{\xi}(A)):=\chi_{\xi}(\gpt_{t}(A))\,, \quad
A\in\ga\,.
$$

As the time evolution and the shift on $\bz^{d}$ are supposed to be
commuting, the random time evolution $\t^{\xi}_{t}$ satisfies the
equivariance condition \eqref{4}, $T_{x}$ being the shift on $X$.

We report the definition of the infinite volume Gibbs condition.
Let $\om\in\cs(\ga)$. The state $\om$ is a Gibbs state at inverse
temperature $\b\neq0$ if
\begin{itemize}
\item[(i)] $\om$ is a modular state,\footnote{This means by
definition, that $\Om_{\om}$ is separating for
$\pi_{\om}(\ga)''$.}
\item[(ii)] for each bounded $\L\subset\bz^{d}$, there exists a
state $\psi_{\L^{c}}\in\cs\big(\ga_{\L^{c}}\big)$ such that
\begin{equation}
\label{gibbs} 
\om^{\b H_{\L}}=\tr_{{}_{\L}}\otimes\psi_{\L^{c}}\,.
\end{equation}
\end{itemize}

Here, $H_{\L}\in\ga$ is the selfadjoint infinitesimal generator of
$\gpt_{t}$ for elements in $\gb_{\L}$ (i.e $\d(A)=i[H_{\L},A]$,
$A\in\gb_{\L}$), $\om^{\b H_{\L}}$ is the perturbation of $\om$ by the 
selfadjoint
element $\b H_{\L}$ (see e.g. \cite{BR2} Theorem 5.4.4), and
finally $\tr_{{}_{\L}}$ is the restriction to $\ga_{\L}$ of the
normalized trace on $\gb$.\footnote{The infinitesimal generator $H_{\L}$
always exists for the situation considered here, taking into
account Proposition 3.3 of \cite{K}, and the construction in
Example 3.2.25 of \cite{BR1}. This can be easily applied to models
with finite range interaction such as that described by the
Hamiltonian \eqref{fham} by taking from the local infinitesimal
generator
$$
H_{\L}:=-\frac{1}{2}\sum_{i\in\L}\sum_{|i-j|=1}J_{i,j}\s(i)\s(j)\,.
$$

Notice that the Gibbs condition 
does depend on the dynamics, as it involves the infinitesimal 
local generators $H_{\L}$.}  The Gibbs condition for the canonical
lifting
\begin{equation*}
\f_{\om}(A):=\int_{G}a_{g}(A)\di g\,, \quad A\in\gb
\end{equation*}
is defined analogously. Denote $\n_{\om}$ the measure on $X$
obtained by restricting $\om$ to $C(X)$, and $\n_{0}$
($\n_{0_{\L}}$) the measure on $X$ ($X_{\L}$) obtained by
restricting the
normalized trace on $\gb$ to $C(X)$ ($C(X_{\L})$).

Some of the main results of \cite{K} are summarized in the
following theorem. We sketch its proof for the reader's
convenience.
\begin{thm}
\label{mkis} 
Let $\om\in\cs(\ga)$. Then the following assertions
are equivalent.
\begin{itemize}
\item[(i)] $\om$ is a Gibbs state at $\b\neq0$,
\item[(ii)] $\f_{\om}\in\cs(\gb)$ is a Gibbs state at $\b\neq0$.
\end{itemize}

If one of the above conditions holds true, then the following
assertions hold true as well.
\begin{itemize}
\item[(a)] $\om$ is a $(\gpt,\b)$--KMS state at $\b\neq0$,
\item[(b)] $\om_{\xi}$ is a $(\t^{\xi},\b)$--KMS
state at $\b\neq0$ almost everywhere w.r.t. the measure
$\n_{\om}$,
\item[(c)] there exists a positive measure $\m$ on $X_{\L^{c}}$ such
that
\begin{equation}
\label{dlrc} 
\om_{\xi}^{\b\chi_{\xi}(H_{\L})}(I)\n_{\om}(\di\xi)
=\n_{0_{\L}}(\di\xi_{\L})\times\m(\di\xi_{\L^{c}})\,.
\end{equation}
\end{itemize}
\end{thm}
\begin{proof}
The equivalence (i)$\iff$(ii) follows from the fact that
$H_{\L}\in\ga$, provided that $\f_{\om}$ is separating if $\om$
is. The last fact is proved in Proposition 4.5 of \cite{K}.

It is well--known that (i)$\Rightarrow$(a). But (a) is equivalent
to (b), see e.g. \cite{K}, or \cite{BF}, Proposition 3.2.

By construction, $\om$ is equivalent to $\om^{\b H_{\L}}$, that is
$$
\om^{\b H_{\L}}=\int^{\oplus}_{X}\big(\om^{\b
H_{\L}}\big)_{\xi}\n_{\om}(\di\xi)\,.
$$

This means that
$$
\n_{\om^{\b H_{\L}}}(\di\xi) =\big(\om^{\b
H_{\L}}\big)_{\xi}(I)\n_{\om}(\di\xi)\,.
$$
Thus, (i)$\Rightarrow$(c) follows by restricting \eqref{gibbs} to
$C(X)$, taking into account that
$$
\big(\om^{\b H_{\L}}\big)_{\xi}=\om_{\xi}^{\b\chi_{\xi}(H_{\L})}
$$
$\n_{\om}$--almost everywhere.
\end{proof}
\begin{rem}
\label{remg}
As the $C^{*}$--algebra $\ca$ is simple, it is well--known that
${\rm(b)}$ of Theorem \ref{mkis} is equivalent to
\begin{itemize}
\item[(b')] $\om_{\xi}$ is a Gibbs state
at $\b\neq0$ $\n_{\om}$--almost everywhere.
\end{itemize}
\end{rem}

Theorem \ref{mkis} leaves open the possibility that the Gibbs
condition could be not equivalent to the KMS condition for a state
$\om\in\cs(\ga)$. This is precisely the question of {\it weak
Gibbsianess}, which was then addressed firstly by
Kishimoto.\footnote{The KMS boundary condition is not sensitive in
classical case. This means that it gives no condition on the
measure $\n_{\om}$. However, it can happen that infinite volume
limits of finite volume Gibbs states could be not Gibbsian on the
skew space describing the joint configurations of spins and
couplings even if $\n_{\om}$ satisfies the DLR condition
\eqref{dlrc}, see below.} Recently, this fact has been checked in \cite{EMSS,
Ku} for classical disordered systems even when the measure
describing the disorder trivially fulfills the
Dobrushin--Lanford--Ruelle (DLR for short, see e.g. \cite{Ru})
condition \eqref{dlrc}. Even if there is no quantum examples
pertaining  to this point, it is expected that the quantum
versions of the models in \cite{EMSS, Ku} exhibit weakly Gibbsian
non jointly Gibbsian states as well.

To end this section, we briefly mention the variational principle
arising from the Gibbs condition for a state $\om\in\ga$. We
suppose that the local Hamiltonians
${\displaystyle\big\{H_{\L}\,\big|\,\L\subset\bz^{d}\,,\L\,\text{bounded}\big\}}$
are associated to a finite range interaction $\Phi$ described by
${\displaystyle\big\{\Phi(\L)\,\big|\,\L\subset\bz^{d}\,,\L\,\text{bounded}\big\}}$,
see e.g. \cite{BR2}. The Gibbs condition for a translation
invariant state $\om\in\cs(\ga)$ is equivalent under the above
conditions to
\begin{equation}
\label{vc1}
p(\Phi,\b)=s(\n_{\om})+\int_{X}p(\chi_{\xi}(\Phi),\b)\n_{\om}(\di\xi)\,,
\end{equation}
\begin{equation}
\label{vc2} p(\chi_{\xi}(\Phi),\b)=s(\om_{\xi})-\b e(\om_{\xi})\,.
\end{equation}

Here, $p(\Phi,\b)$ ($p(\chi_{\xi}(\Phi),\b)$) is the pressure
associated with the interaction $\Phi$ ($\chi_{\xi}(\Phi)$
described by the continuous field of potentials\\ 
${\displaystyle\big\{
\chi_{\xi}(\Phi(\L))\,\big|\,\L\subset\bz^{b}\,,\L\,\text{bounded}\big\}}$),
$s(\n_{\om})$ is the mean entropy of the translationally invariant
measure $\n_{\om}$, $s(\om_{\xi})$ and $e(\om_{\xi})$ the mean
entropy and the mean energy of the state $\om_{\xi}$,
respectively.\footnote{See e.g. \cite{BR2} for the definition of
the pressure. Furthermore, $s(\om_{\xi})$ and $e(\om_{\xi})$
should be defined with a little bit of care as the $\om_{\xi}$ are not
translationally invariant, see Prosition 7.5 and Prosition 7.6 of
\cite{K}.} The functions $p(\chi_{\xi}(\Phi),\b)$, $s(\om_{\xi})$,
and $e(\om_{\xi})$ are measurable essentially bounded functions
defined $\n_{\om}$--almost everywhere, so they depend also on the 
state $\om$ under consideration. Thus, \eqref{vc2} holds true 
almost surely w.r.t the measure $\n_{\om}$. Notice that, if $\n_{\om}$ is
ergodic w.r.t. the spatial translations, they are constant 
$\n_{\om}$--almost
everywhere. Equation \eqref{vc2} is nothing but the variational 
principle associated to the Gibbs condition for the states 
$\om_{\xi}$, which holds true $\n_{\om}$--almost everywhere, taking 
into account Remark \ref{remg}.

We refer the reader to \cite{A1, K} for the meaning of
\eqref{vc1} and  \eqref{vc2} as variational principles, and for
further details.

\end{document}